\documentclass[12pt]{article}

\usepackage{amssymb}
\usepackage{amsmath,amsthm}
\usepackage[dvips]{graphicx}
\usepackage{wrapfig}
\usepackage{bm}

\newtheorem{theorem}{Theorem}[section]

\newtheorem{proposition}[theorem]{Proposition}

\theoremstyle{definition}

\theoremstyle{remark}

\title{The reductivity of spherical curves Part I\hspace{-1.5pt}I: 4-gons}
\author{
Yui Onoda \thanks{Department of School Education, Akita University, 1-1 Akita-Gakuenmachi, Akita 010-8552, Japan. }
\and 
Ayaka Shimizu \thanks{Department of Mathematics, National Institute of Technology, Gunma College, 580 Toriba-cho, Maebashi-shi, Gunma 371-8530, Japan. Email: shimizu@nat.gunma-ct.ac.jp }
}
\date{\today}

\begin{document}

\maketitle

\begin{abstract}
The reductivity of a spherical curve represents how reduced it is. 
It is unknown if there exists a spherical curve whose reductivity is four. 
In this paper we give an unavoidable set for spherical curves with reductivity four, that is, we give a set of parts of spherical curves such that every spherical curve with reductivity four has at least one of the parts,  by considering 4-gons. 
\end{abstract}

\section{Introduction}

A {\it spherical curve}, or a {\it knot projection}, is a closed curve on the 2-sphere, where each crossing is a double point and crosses transversely. 
In this paper we consider spherical curves with at least one crossing, and we consider spherical curves up to reflection, that is, we assume the spherical curve which is obtained by turning over a spherical curve $P$ and $P$ are the same.  
We call each part of the sphere bounded by a spherical curve a {\it region}. 
In particular, we call a region which has $n$ edges on the boundary an {\it $n$-gon}. 
A crossing of a spherical curve is {\it reducible} if there are just three regions around the crossing. 
A spherical curve is said to be {\it reducible} if it has a reducible crossing. 
If not, then it is said to be {\it reduced}. 
The reductivity $r(P)$, defined by the second author in \cite{shimizu}, of a spherical curve $P$ represents how far $P$ is from a reducible spherical curve. 
The precise definition of reductivity will be given in Section 2. 
It is shown in \cite{shimizu} that the reductivity is four or less for every spherical curve although we have not found a spherical curve with reductivity four. 
Our next step is to prove that the reductivity is three or less for every spherical curve, or to find a spherical curve whose reductivity is four. \\

It is shown in \cite{AST} that every reduced spherical curve has a bigon (2-gon) or trigon (3-gon). 
This is because we have the formula 
\begin{align}
2 C_2 + C_3 = 8+ C_5 + 2 C_6 +3 C_7 +4 C_8 + \dots 
\label{238formula}
\end{align}
for every reduced spherical curve, where $C_n$ is the number of $n$-gons. 
\begin{figure}[h]
\begin{center}
\includegraphics[width=90mm]{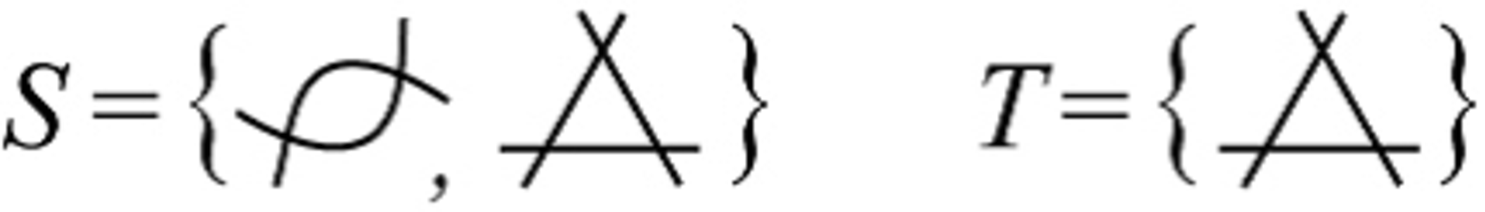}
\caption{$S$ is an unavoidable set for reduced spherical curves. The set $T$ is an unavoidable set for spherical curves with reductivity three or four. }
\label{st}
\end{center}
\end{figure}
In other words, the set $S$ in Fig.\ref{st} is an {\it unavoidable set} for reduced spherical curves. 
As mentioned in \cite{shimizu}, the reductivity of a spherical curve $P$ is two or less if $P$ has a bigon. 
Therefore we can say that the set $T$ in Fig.\ref{st} consisting of a trigon is an unavoidable set for spherical curves with reductivity three or four. 
Moreover, we can say that every spherical curve with reductivity three or four has at least eight trigons because of the formula (\ref{238formula}). \\

\begin{figure}[h]
\begin{center}
\includegraphics[width=100mm]{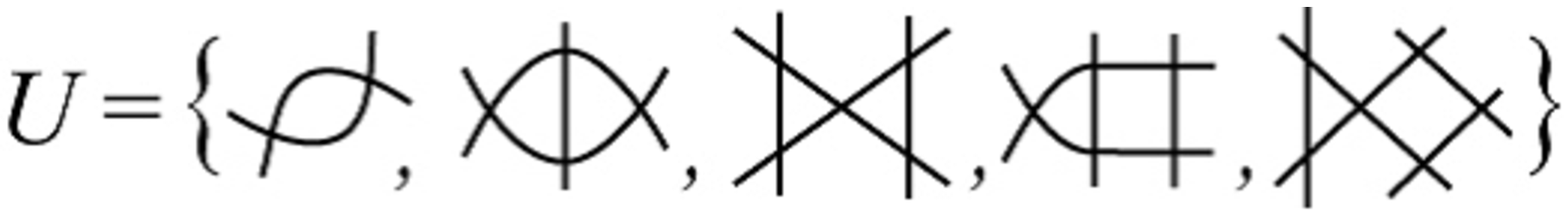}
\caption{$U$ is an unavoidable set for reduced spherical curves. }
\label{s1}
\end{center}
\end{figure}

Using the ``discharging method'' with the formula (\ref{238formula}) (inspired by the story of the four color theorem), it is shown in \cite{shimizu} that the set $U$ shown in Fig.\ref{s1} is also an unavoidable set for reduced spherical curves, that is, every reduced spherical curve has at least one part in $U$. 
In \cite{shimizu}, trigons of a spherical curve are classified into four types A, B, C and D (see Fig.\ref{2and3-di}), and in Section 3 of this paper, 4-gons of a spherical curve are classified into 13 types $1a$, $1b$, $\dots$, $4b$ (see Fig.\ref{4-di}) with respect to the outer connections. 
In this paper we show the following: 

\phantom{x}
\begin{theorem}
If there exists a spherical curve $P$ with reductivity four, $P$ has at least one of the 21 parts depicted in Fig.\ref{unavoidable-r4}. 
\begin{figure}[h]
\begin{center}
\includegraphics[width=120mm]{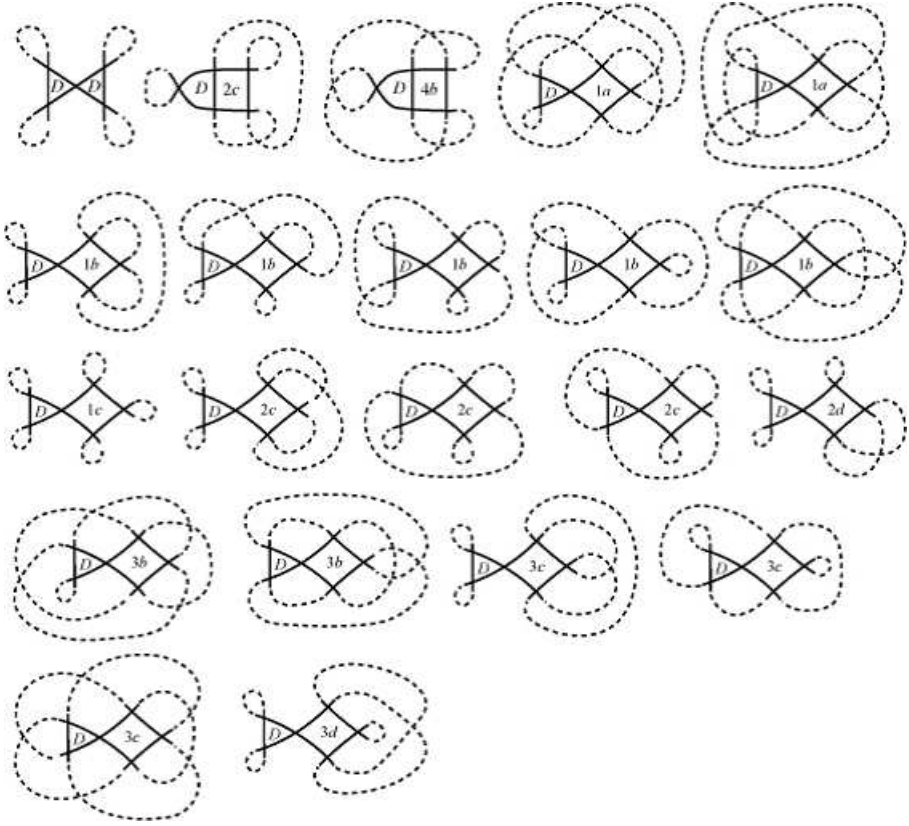}
\caption{An unavoidable set for spherical curves with reductivity four.}
\label{unavoidable-r4}
\end{center}
\end{figure}
\label{thm-r4}
\end{theorem}
\phantom{x}

\noindent Note that dotted arcs in the figures represent the outer connection in this paper. 
The theorem above implies that the set consisting of the 21 parts in Fig.\ref{unavoidable-r4} is an unavoidable set for spherical curves with reductivity four. 

The rest of this paper is organized as follows: 
In Section 2, we review the reductivities. 
In Section 3, we classify 4-gons and prove Theorem \ref{thm-r4}.

\section{Reductivity}

\begin{figure}[h]
\begin{center}
\includegraphics[width=80mm]{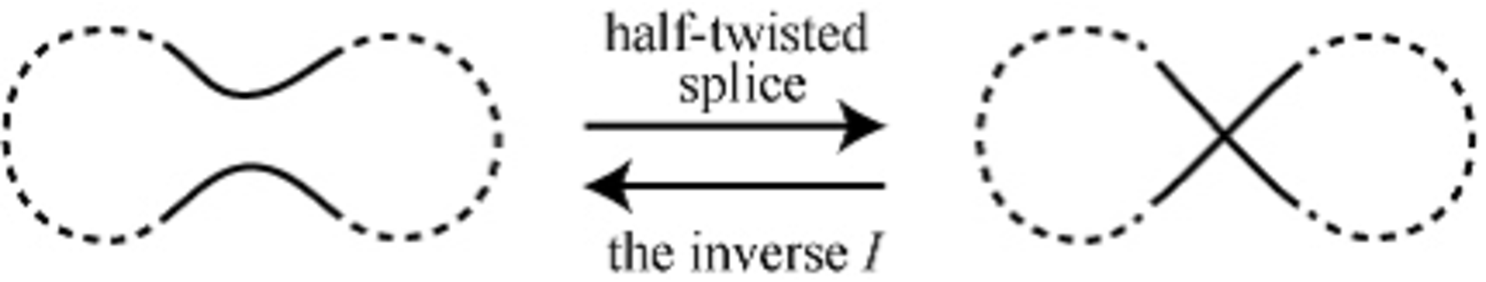}
\caption{Half-twisted splice and its inverse $I$.}
\label{halftwisted}
\end{center}
\end{figure}

In this section we review the reductivity. 
A {\it half-twisted splice} is a local transformation on a spherical curve illustrated in Fig.\ref{halftwisted} (\cite{calvo, ito-shimizu}). 
We call $I$ the inverse of a half-twisted splice. 
The {\it reductivity}, denoted by $r(P)$, of a spherical curve $P$ is the minimal number of $I$ required to obtain a reducible spherical curve from $P$. 

\begin{figure}[h]
\begin{center}
\includegraphics[width=120mm]{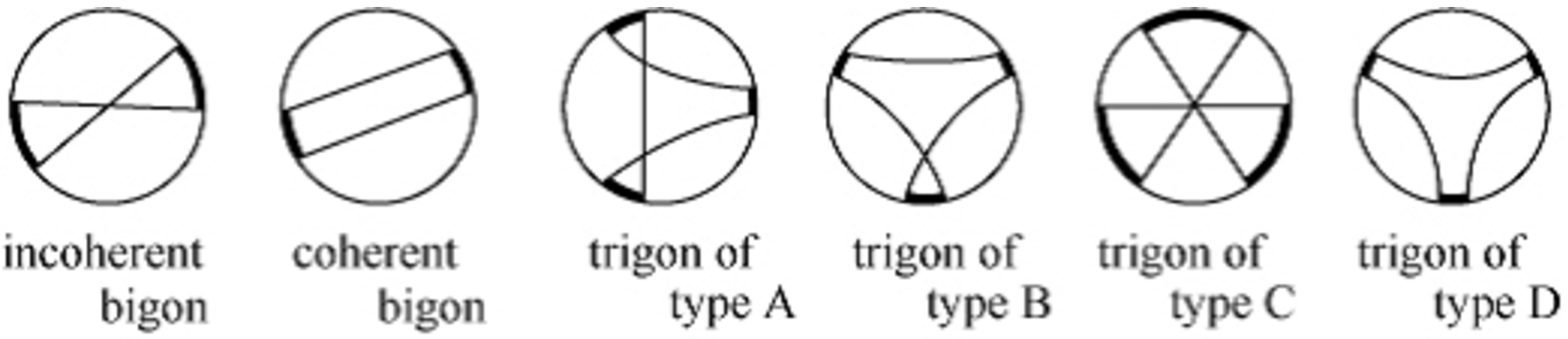}
\caption{Bigons and trigons on chord diagrams. There are no segments on the thick arcs. }
\label{2and3-co}
\end{center}
\end{figure}

Considering the outer connections, bigons and trigons are divided into the two and four types as illustrated in Fig.\ref{2and3-di}. 
A {\it chord diagram} of a spherical curve is the preimage of the spherical curve with two points corresponding to the same crossing connected with a segment when we consider the spherical curve as an immersion of a circle to the sphere. 
The bigons and trigons on chord diagrams are illustrated in Fig.\ref{2and3-co}. 
The following proposition is shown in \cite{shimizu}: 
\begin{figure}[h]
\begin{center}
\includegraphics[width=120mm]{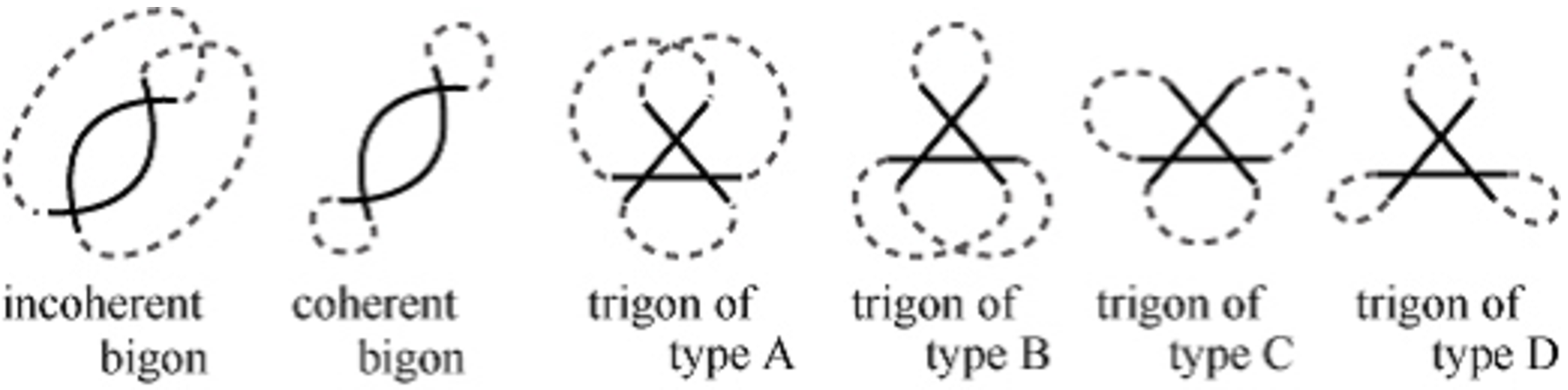}
\caption{Bigons and trigons.}
\label{2and3-di}
\end{center}
\end{figure}

\phantom{x}
\begin{proposition}[\cite{shimizu}]
Let $P$ be a spherical curve. 
\begin{itemize}
\item If $P$ has an incoherent bigon, $r(P) \le 1$. 
\item If $P$ has a coherent bigon, $r(P) \le 2$. 
\item If $P$ has a trigon of type $A$, $r(P) \le 2$. 
\item If $P$ has a trigon of type $B$, $r(P) \le 3$. 
\item If $P$ has a trigon of type $C$, $r(P) \le 3$. 
\end{itemize}
\label{bitri}
\end{proposition}
\phantom{x}

\begin{figure}[h]
\begin{center}
\includegraphics[width=50mm]{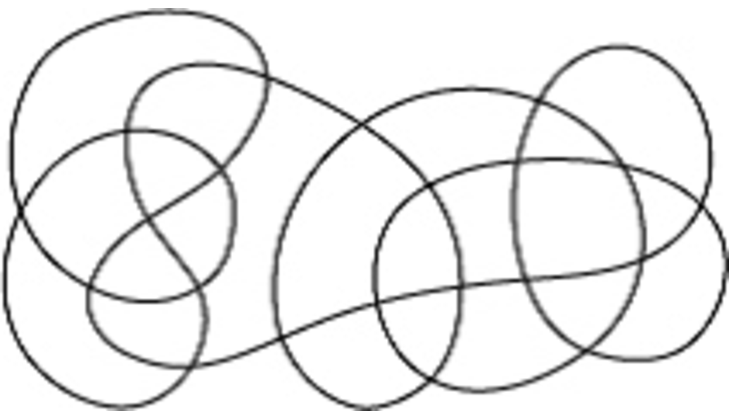}
\caption{A spherical curve without bigons whose reductivity is one.}
\label{no-bi}
\end{center}
\end{figure}

\begin{figure}[h]
\begin{center}
\includegraphics[width=50mm]{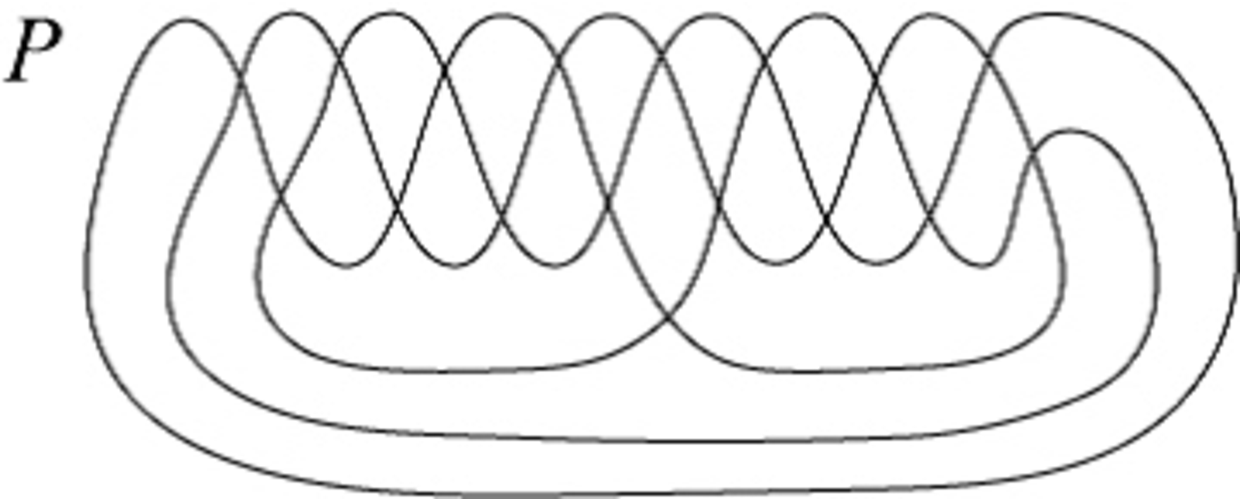}
\caption{The reductivity of $P$ is three.}
\label{r3ex}
\end{center}
\end{figure}

\noindent Remark that the converse of Proposition \ref{bitri} does not hold. 
For example, the spherical curve in Fig.\ref{no-bi} has no bigons, but the reductivity is one. 
In \cite{ito-takimura}, the definitions of spherical curves with reductivity one and two are rephrased diagramatically (see Fact 1 and Theorem 1 in \cite{ito-takimura}) and that is useful to determine reductivity. 
For example, we can say that the spherical curve $P$ depicted in Fig.\ref{r3ex} has the reductivity three because $r(P)$ is not two or less, and $r(P) \le 3$ because $P$ has trigons of type B and C. 
We also have the following proposition: 

\phantom{x}
\begin{proposition}
The standard projection of a $(3, 3n+1)$-torus knot has reductivity three ($n=1,2,3, \dots $).
\label{3n+1}
\end{proposition}
\phantom{x}

\noindent Note that all the trigons of the standard projection of a $(3, 3n+1)$-torus knot are type B. 
From the proposition above, we can say that there are infinitely many spherical curves with reductivity three. 
Similarly, we also have the following. 

\phantom{x}
\begin{proposition}
The standard projection of a $(3, 3n-1)$ and $(4, 4n-1)$-torus knot has the reductivity two. 
The standard projection of a $(4, 4n+1)$-torus knot has the reductivity three ($n=1,2,3, \dots $). 
\label{3n-1}
\end{proposition}
\phantom{x}

\noindent Note that the standard projection of a $(3, 3n-1)$ and $(4, 4n-1)$-torus knot has a trigon of type A.

\begin{figure}[h]
\begin{center}
\includegraphics[width=80mm]{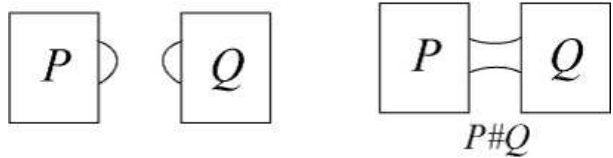}
\caption{Connected sum.}
\label{conn}
\end{center}
\end{figure}

For spherical curves $P$ and $Q$, we call a spherical curve obtained from $P$ and $Q$ by connecting their arcs as depicted in Fig.\ref{conn} a {\it connected sum} of $P$ and $Q$, and denote it by $P \sharp Q$. 
We say a spherical curve is {\it prime} if it is not any connected sum of two (nontrivial) spherical curves. 
Since one part can be regarded as a small 1-tangle for the other one, we have the following inequality: 
$$
r(P \sharp Q ) \le \text{min} \{ r(P), r(Q) \}.
$$

\noindent In \cite{oregon}, the chord diagrams for all the prime spherical curves with up to ten crossings are listed. 
We have the following: 

\phantom{x}
\begin{proposition}
All the spherical curves with 10 or less crossings except the standard projection of a $(3,4)$-torus knot have the reductivity two or less. 
\end{proposition}
\phantom{x}

\begin{proof}
All the prime spherical curves with ten or less crossings except the standard projections of the $(3,4), (4,3)$ and $(3,5)$-torus knot (see Fig.\ref{34-torus} and \ref{43-35-torus}) have a bigon. 
By Propositions \ref{bitri}, \ref{3n+1} and \ref{3n-1}, the claim holds for prime spherical curves. 
For non-prime spherical curves, since they are connected sums of spherical curves up to nine crossings, the claim holds. 
Note that even if it is a connected sum of the standard projection $P$ of the $(3,4)$-torus knot and some spherical curves, the reductivity is zero because it is a connected sum of $P$ (which has eight crossings) and spherical curves with one or two crossings. 
\begin{figure}[h]
\begin{center}
\includegraphics[width=40mm]{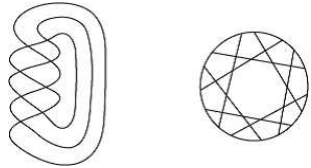}
\caption{The standard projection of the $(3,4)$-torus knot.}
\label{34-torus}
\end{center}
\end{figure}

\begin{figure}[h]
\begin{center}
\includegraphics[width=80mm]{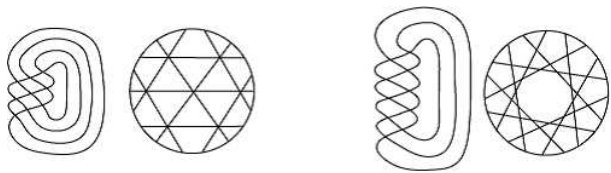}
\caption{The standard projections of the $(4,3)$ and $(3,5)$-torus knot.}
\label{43-35-torus}
\end{center}
\end{figure}

\end{proof}

\section{4-gons}

\begin{figure}[h]
\begin{center}
\includegraphics[width=80mm]{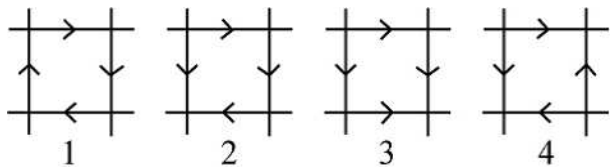}
\caption{4-gons with relative orientations.}
\label{4-ori}
\end{center}
\end{figure}

\begin{figure}[h]
\begin{center}
\includegraphics[width=120mm]{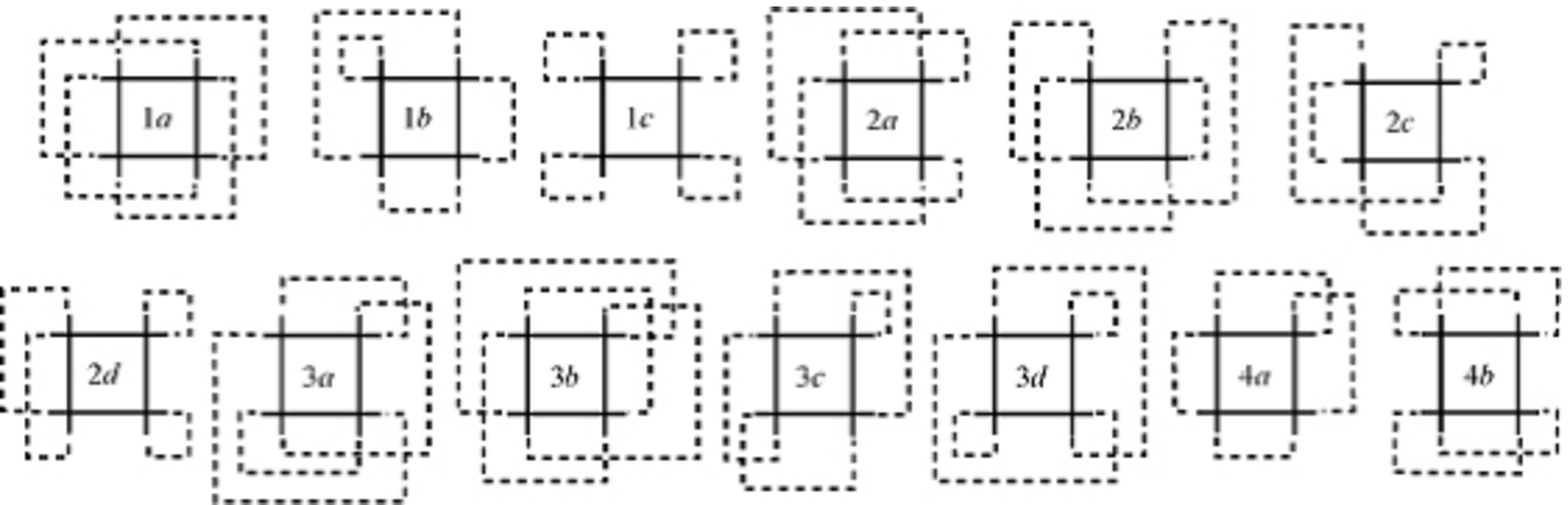}
\caption{The 13 types of 4-gons.}
\label{4-di}
\end{center}
\end{figure}

\begin{figure}[h]
\begin{center}
\includegraphics[width=120mm]{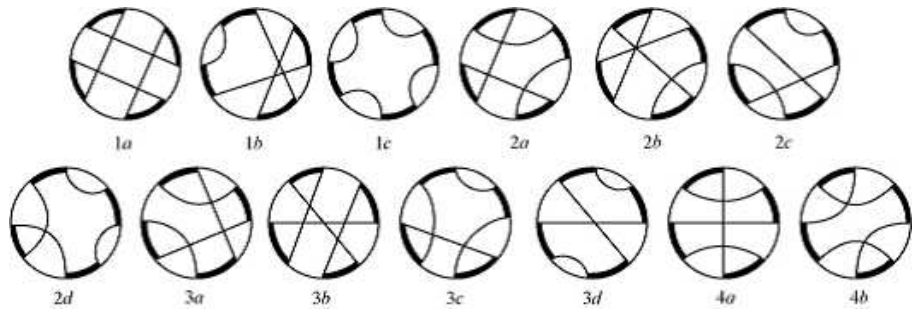}
\caption{The 4-gons on chord diagrams. There are no segments on the thick arcs. }
\label{4-co}
\end{center}
\end{figure}

In this section we classify 4-gons and prove Theorem \ref{thm-r4}. 
We first classify 4-gons into type 1 to 4 with respect to the relative orientations of the boundary as illustrated in Fig.\ref{4-ori}. 
Then we classify them with respect to the outside connection and name them $1a, 1b, 1c, \dots $ and $4b$ as shown in Fig.\ref{4-di}. 
See Fig.\ref{4-co} for chord diagrams. 

\noindent We have the following: 

\phantom{x}
\begin{proposition}
If a spherical curve $P$ has a 4-gon of type $2a, 2b, 3a$ or $4a$, then $r(P) \le 3$. 
\label{prop3-4gon}
\end{proposition}
\phantom{x}

\begin{proof}

The 4-gons of type $2a, 2b, 3a$ and $4a$ are illustrated in Fig.\ref{3d-4-di}. 
\begin{figure}[h]
\begin{center}
\includegraphics[width=80mm]{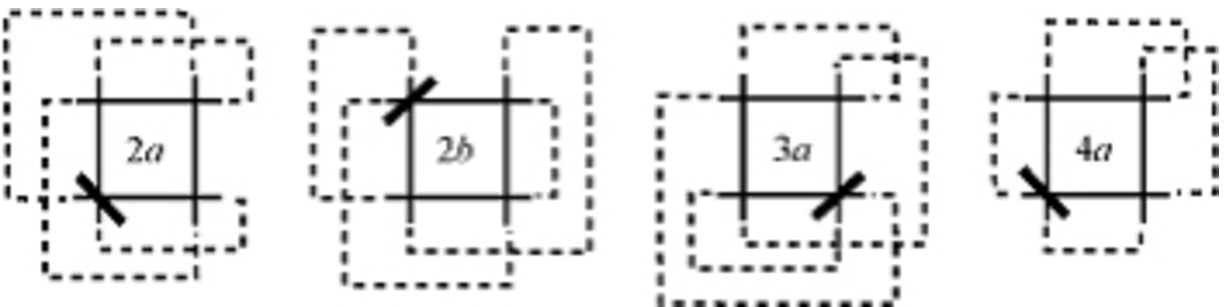}
\caption{The 4-gons we can obtain a trigon of type A by an $I$.}
\label{3d-4-di}
\end{center}
\end{figure}
By applying $I$ at the marked crossing, we obtain a trigon of type A from each 4-gons. 
\end{proof}
\phantom{x}

\noindent Now we prove Theorem \ref{thm-r4}. 

\phantom{x}
\noindent {\it Proof of Theorem \ref{thm-r4}.} \ 
If there exists a spherical curve $P$ with $r(P)=4$, $P$ has at least one parts of $U$ in Fig.\ref{s1} because $P$ is a reduced spherical curve. 
Here, $P$ can not have any bigon, trigon of type A, B or C, or 4-gon of type $2a, 2b, 3a$ or $4a$ by Propositions \ref{bitri} and \ref{prop3-4gon}. 
Then $P$ can not have the first part of $U$ because it is a bigon. 
Also, $P$ can not have the second part of $U$ because even if one of the two trigons is of type D, then another one should be a trigon of type A or B because the boundary is incoherent. 
Hence it is sufficient to consider the rest three parts of $U$. 

\begin{figure}[h]
\begin{center}
\includegraphics[width=30mm]{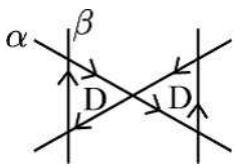}
\caption{The third part of $U$.}
\label{3ofu}
\end{center}
\end{figure}

\phantom{x}
If $P$ has the third part of $U$, the both two trigons should be of type D, and we have just one suitable outer connection; 
For example, in Fig.\ref{3ofu}, $\beta$ can be connected outside with only $\alpha$ considering the outer connections of the both two trigons. 
Thus we obtain the 1st part in Fig.\ref{unavoidable-r4}. 

\begin{figure}[h]
\begin{center}
\includegraphics[width=30mm]{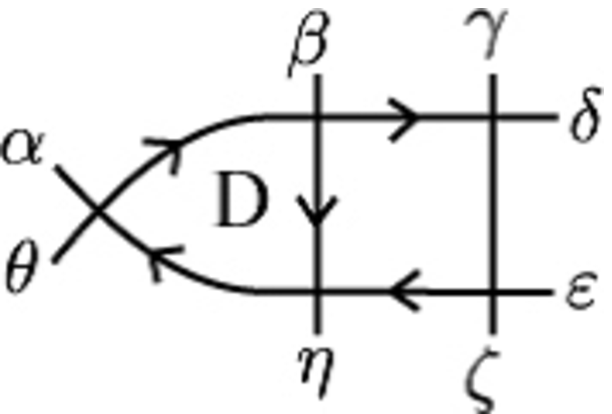}
\caption{The fourth part of $U$.}
\label{4ofu}
\end{center}
\end{figure}

\phantom{x}
If $P$ has the fourth part of $U$, the trigon should be of type D, and then the orientation of the 4-gon should be of type 2 or 4 (see Fig.\ref{4ofu}). 
More precisely, the 4-gon should be of type $2c$, $2d$ or $4b$ by Proposition \ref{prop3-4gon}. 

\begin{figure}[h]
\begin{center}
\includegraphics[width=70mm]{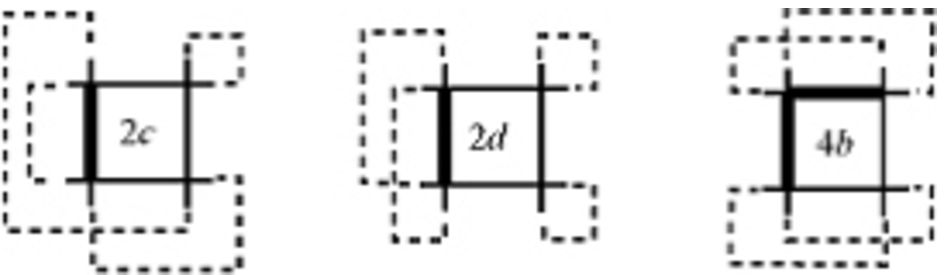}
\caption{The 4-gon of the fourth part should be of type $2c$, $2d$ or $4b$, and a thick edge is to be shared with the trigon.}
\label{44edge}
\end{center}
\end{figure}

\noindent $\bullet$ For the case of 4-gon of type $2c$, the thick edge of the 4-gon of type $2c$ in Fig.\ref{44edge} is to be shared with the trigon because of the orientation. 
Now we consider the outer connection. 
See Fig.\ref{4ofu}. 
Considering the outer connection of the trigon of type D and the 4-gon of type $2c$, $\delta$ and $\gamma$ should be connected outside. 
Then $\zeta$ and $\beta$, $\eta$ and $\varepsilon$ should be connected outside. 
Thus we have just one suitable outer connection of a trigon of type D and a 4-gon of type $2c$, which is the 2nd part in Fig.\ref{unavoidable-r4}. 

\noindent $\bullet$ Next we consider the case of 4-gon of type $2d$. 
The thick edge of the 4-gon of type $2d$ in Fig.\ref{44edge} is to be shared with the trigon. 
Considering the outer connection of the trigon of type $D$, $\alpha$ should be connected outside with $\theta$ or $\gamma$. 
On the other hand, considering the outer connection of the 4-gon of type $2d$, $\alpha$ should be connected outside with $\beta$. 
Hence there are no suitable connection in this case. 

\noindent $\bullet$ Next we consider the case of 4-gon of type $4b$. 
From a viewpoint of the orientation, each edge of a 4-gon of type $4b$ may be shared with the trigon, and it is sufficient to consider the two thick edges in Fig.\ref{44edge} by symmetry. 

First, we consider the case the left-side edge of the 4-gon of type $4b$ in Fig.\ref{44edge} is shared with the trigon. 
Considering the outer connection of the trigon of type D and the 4-gon of type $4b$, $\alpha$ and $\zeta$ should be connected outside. 
Then $\gamma$ and $\theta$, $\delta$ and $\beta$, and $\eta$ and $\varepsilon$ should be connected outside. 
Thus, in this case we obtain just one suitable connection, which is the 3rd part in Fig.\ref{unavoidable-r4}. 

Next we consider the case the upper-side edge of the 4-gon of type $4b$ in Fig.\ref{44edge} is shared with the trigon. 
Considering the outside connection of the trigon of type D, $\alpha$ should be connected outside with $\theta$ or $\zeta$. 
On the other hand, considering the outer connection of the 4-gon of type $4b$, $\alpha$ should be connected outside with $\beta$. 
Hence there are no suitable connection in this case.

\begin{figure}[h]
\begin{center}
\includegraphics[width=30mm]{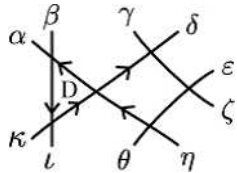}
\caption{The fifth part of $U$.}
\label{5ofu}
\end{center}
\end{figure}

\phantom{x}
If $P$ has the fifth part of $U$, the trigon also should be of type D, and then the orientation of the 4-gon should be of type 1, 2 or 3 (see Fig.\ref{5ofu}). 
More precisely, the 4-gon should be of type $1a, 1b, 1c, 2c, 2d, 3b, 3c$ or $3d$ by Proposition \ref{prop3-4gon}. 

\begin{figure}[h]
\begin{center}
\includegraphics[width=100mm]{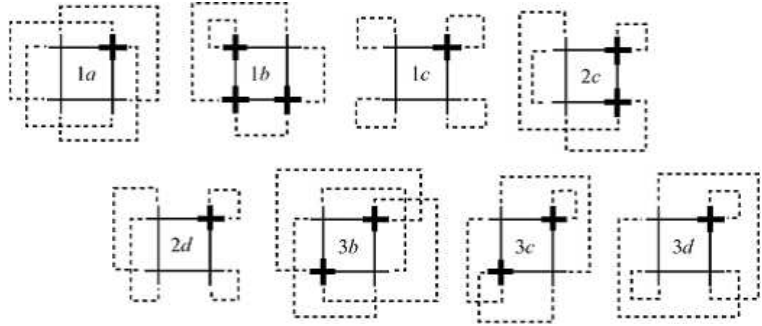}
\caption{The 4-gon of the fifth part should be of type $1a, 1b, 1c, 2c, 2d, 3b, 3c$ or $3d$, and a thick crossing is to be shared with the trigon.}
\label{55edge}
\end{center}
\end{figure}

\noindent $\bullet$ We first consider the case of 4-gon of type $1a$. 
By the orientation and the symmetry, it is sufficient to check the case the thick crossing of the 4-gon of type $1a$ in Fig.\ref{55edge} is shared with the trigon. 
See Fig.\ref{5ofu}. 
Considering the outer connection of the trigon and 4-gon, $\delta$ should be connected outside with $\eta$. 
Then, $\alpha$ should be connected outside with $\beta$ or $\varepsilon$. 

For the case $\alpha$ is connected with $\beta$, $\iota$ is to be connected outside with $\varepsilon$ and then $\theta$ and $\gamma$, and $\zeta$ and $\kappa$ should be connected outside. 
Thus we obtain a suitable connection which is the 4th part in Fig.\ref{unavoidable-r4}. 

For the case $\alpha$ is connected outside with $\varepsilon$, then $\theta$ should be connected outside with $\beta$ or $\gamma$. 

If $\theta$ and $\beta$ are connected outside, then $\iota$ and $\gamma$, and $\zeta$ and $\kappa$ should be connected outside. 
Thus we obtain a suitable connection, which is the 5th part in Fig.\ref{unavoidable-r4}. 

Next we consider the rest case on the 4-gon of type $1a$ that $\delta$ is connected outside with $\eta$, $\alpha$ is connected outside with $\varepsilon$, and $\theta$ is connected outside with $\gamma$. 
In this case, $\zeta$ should be connected outside with $\beta$, and $\iota$ should be conected outside with $\kappa$. 
Thus we obtain a suitable connection, which is also the 4th part in Fig.\ref{unavoidable-r4}.

\noindent $\bullet$ Next is the case of 4-gon of type $1b$. 
The three thick crossings of the 4-gon of type $1b$ in Fig.\ref{55edge} can be shared with the trigon. 

First we consider the case the left-upper crossing is shared with the trigon. 
Considering the outer connection, $\delta$ should be connected outside with $\varepsilon$. 
Then $\theta$ and $\gamma$, $\zeta$ and $\eta$, $\alpha$ and $\beta$, and $\iota$ and $\kappa$ should be connected outside. 
Thus we obtain a suitable connection, which is the 6th one. 

Next we consider the case the left-lower crossing is shared with the trigon. 
In this case, $\delta$ should be connected outside with $\varepsilon$, and $\theta$ should be connected outside with $\eta$. 
Then, $\alpha$ is to be connected outside with $\gamma$ or $\beta$. 

For the case $\alpha$ and $\gamma$ are connected outside, $\zeta$ and $\beta$, and $\iota$ and $\kappa$ should be connected outside. 
Thus we obtain the 7th one. 

For the case $\alpha$ and $\beta$ are connected outside, $\iota$ and $\gamma$, and $\zeta$ and $\kappa$ should be connected outside. 
Thus we obtain the 8th one. 

Finally we consider the case the right-lower crossing of the 4-gon of type $1b$ is shared with the trigon. 
First, $\delta$ should be connected outside with $\eta$. 
Then $\alpha$ is to be connected outside with $\gamma$ or $\beta$. 

For the case $\alpha$ and $\gamma$ are connected outside, then $\zeta$ is to be connected outside with $\varepsilon$ or $\beta$. 

For the case $\zeta$ and $\varepsilon$ are connected outside, then $\theta$ and $\beta$, and $\iota$ and $\kappa$ should be connected outside. 
Thus we obtain the 9th one. 

Next we consider the case $\zeta$ and $\beta$ are connected outside. 
In this case, $\iota$ and $\varepsilon$, and $\theta$ and $\kappa$ should be connected outside. 
Thus we obtain the 10th one. 

Finally, we consider the case that $\delta$ and $\eta$, and $\alpha$ and $\beta$ are connected outside. 
In this case, $\iota$ and $\gamma$, $\zeta$ and $\varepsilon$, and $\theta$ and $\kappa$ should be connected outside. 
Thus we obtain a suitable connection, which is also the 9th one. 

\noindent $\bullet$ Next, we consider the case of 4-gon of type $1c$. 
By symmetry, it is sufficient to consider the case the thick crossing in Fig.\ref{55edge} is shared with the trigon. 
In this case, $\delta$ and $\gamma$, $\zeta$ and $\varepsilon$, $\theta$ and $\eta$, $\alpha$ and $\beta$, and $\iota$ and $\kappa$ should be connected outside. 
Hence we have only one suitable connection, which is the 11th one. 

\noindent $\bullet$ Next we consider the case of 4-gon of type $2c$. 
There are two thick crossings which can be shared with the trigon. 

First let us consider the case the right-upper crossing is shared with the trigon. 
In this case, $\delta$ and $\theta$, $\varepsilon$ and $\gamma$, $\zeta$ and $\eta$, $\alpha$ and $\beta$, and $\iota$ and $\kappa$ should be connected outside. 
We obtain the 12th one. 

Next we consider the case the right-lower crossing is shared with the trigon. 
In this case, $\delta$ and $\varepsilon$, and $\theta$ and $\eta$ should be connected outside. 
Then $\alpha$ should be connected outside with $\zeta$ or $\beta$. 

For the case $\alpha$ and $\zeta$ are connected outside, then $\gamma$ and $\beta$, and $\iota$ and $\kappa$ should be connected outside. 
Thus we obtain the 13th one. 

For the case $\alpha$ and $\beta$ are connected outside, $\iota$ and $\zeta$, and $\gamma$ and $\kappa$ should be connected outside. 
Thus we obtain the 14th one. 

\noindent $\bullet$ Next we consider the case of 4-gon of type $2d$. 
By symmetry, it is sufficient to consider the case the thick crossing in Fig.\ref{55edge} is shared with the trigon. 
In this case, $\delta$ and $\gamma$, $\zeta$ and $\theta$, $\varepsilon$ and $\eta$, $\alpha$ and $\beta$, and $\iota$ and $\kappa$ should be connected outside. 
Thus we obtain the 15th one. 

\noindent $\bullet$ Next we consider the case of 4-gon of type $3b$. 
The two thick crossings in Fig.\ref{55edge} can be shared with the trigon. 

First let us consider the case the right-upper crossing is shared with the trigon. 
In this case, $\delta$ should be connected outside with $\eta$, and then $\alpha$ should be connected outside with $\zeta$ or $\beta$. 

For the case $\alpha$ and $\zeta$ are connected outside, $\gamma$ should be connected outside to $\theta$ or $\beta$. 

When $\gamma$ and $\theta$ are connected outside, then $\varepsilon$ and $\beta$, and $\iota$ and $\kappa$ should be connected outside. 
Thus we obtain the 16th one. 

When $\gamma$ and $\beta$ are connected outside, then $\iota$ and $\theta$, and $\varepsilon$ and $\kappa$ should be connected outside. 
Thus we obtain the 17th one. 

For the case $\alpha$ and $\beta$ are connected outside, then $\iota$ and $\zeta$, $\gamma$ and $\theta$, and $\varepsilon$ and $\kappa$ should be connected outside. 
Again, we obtain the 16th one. 

Next, we consider the case that the left-lower crossing is shared with the trigon. 
In this case, $\eta$ should be connected outside with $\delta$. 
Then $\kappa$ should be connected outside with $\varepsilon$ or $\iota$. 

For the case that $\kappa$ and $\varepsilon$ are connected outside, then $\theta$ should be connected outside to $\gamma$ or $\iota$. 

When $\theta$ and $\gamma$ are connected, $\zeta$ and $\iota$, and $\beta$ and $\alpha$ should be connected outside. 
We obtain the 16th one again. 

When $\theta$ and $\iota$ are connected, $\beta$ and $\gamma$, and $\zeta$ and $\alpha$ should be connected outside. 
We obtain the 17th one again. 

Finally we consider the case that $\kappa$ and $\iota$ are connected. 
In this case, $\beta$ and $\varepsilon$, $\theta$ and $\gamma$, and $\zeta$ and $\alpha$ should be connected outside. 
Then we obtain the 16th one again. 

\noindent $\bullet$ Next is the case of 4-gon of type $3c$. 
The two thick crossings in Fig.\ref{55edge} can be shared with the trigon. 

For the case that the right-upper crossing is shared with the trigon, $\delta$ and $\zeta$, $\gamma$ and $\theta$, $\varepsilon$ and $\eta$, $\alpha$ and $\beta$, and $\iota$ and $\kappa$ should be connected outside. 
This is the 18th one. 

For the case that the left-lower crossing is shared with the trigon, $\eta$ should be connected outside with $\delta$. 
Then $\kappa$ should be connected outside with $\gamma$ or $\iota$. 

For the case $\kappa$ is connected with $\gamma$, $\zeta$ should be connected outside with $\varepsilon$ or $\iota$. 

When $\zeta$ and $\varepsilon$ are connected, then $\theta$ and $\iota$, and $\beta$ and $\alpha$ should be connected. 
Then we obtain the 19th one. 

When $\zeta$ and $\iota$ are connected, then $\beta$ and $\varepsilon$, and $\theta$ and $\alpha$ should be connected. 
Then we obtain the 20th one. 

Next we consider the case $\kappa$ and $\iota$ are connected outside. 
In this case $\beta$ and $\gamma$, $\zeta$ and $\varepsilon$, and $\theta$ and $\alpha$ should be connected outside. 
We obtain the 19th one again. 

\noindent $\bullet$ Finally, we consider the case of 4-gon of type $3d$. 
By symmetry, it is sufficient to consider the case of the thick crossing in Fig.\ref{55edge} is shared with the trigon. 
In this case, $\delta$ and $\theta$, $\varepsilon$ and $\zeta$, $\gamma$ and $\eta$, $\alpha$ and $\beta$, and $\iota$ and $\kappa$ should be connected outside. 
Thus, we obtain the 21st one of $U$.

\hfill$\square$

\phantom{x}

\section*{Acknowledgments}
The authors would like to thank the referee and editor for helpful suggestions. 
The second author was partially supported by the Sumitomo Foundation.


\begin{thebibliography}{9}
\bibitem{AST} \textsc{C.~C.~Adams, R.~Shinjo, K.~Tanaka}, Complementary regions of knot and link diagrams, {\it Ann. Comb.} {\bf 15} (2011), 549--563. 
\bibitem{calvo} \textsc{J.~A.~Calvo}, Knot enumeration through flypes and twisted splices, J. Knot Theory Ramifications {\bf 6} (1997), 785--798. 
\bibitem{oregon} \textsc{M.~Chmutov, T.~Hulse, A.~Lum, P.~Rowell}, Plane and spherical curves: an investigation of their invariants, Proceedings of the Research Experiences for Undergraduates Program in Mathematics, Oregon State University, 2006, 1--92. 
\bibitem{ito-shimizu} \textsc{N.~Ito and A.~Shimizu}, The half-twisted splice operation on reduced knot projections, J. Knot Theory Ramifications {\bf 21}, 1250112 (2012) [10 pages]. 
\bibitem{ito-takimura} \textsc{N.~Ito and Y.~Takimura}, Knot projections with reductivity two, Topology Appl. {\bf 193} (2015), 290--301. 
\bibitem{shimizu} \textsc{A.~Shimizu}, The reductivity of spherical curves, Topology Appl. {\bf 196} (2015), 860--867. 
\end{thebibliography}
\end{document}